# Modeling localized failure in geomaterials by large-deformation-plasticity periporomechanics

Xiaoyu Song* and Hossein Pashazad

Engineering School of Sustainable Infrastructure and Environment, University of Florida, Gainesville, Florida, USA

*xysong@ufl.edu

**Abstract.** Large-deformation localized failure in geomaterials plays a crucial role in geohazards engineering. This article investigates shear bands and retrogressive failure of geomaterials through a recently formulated large-deformation-plasticity periporomechanics (PPM) paradigm. Periporomechanics is a nonlocal reformulation of classical poromechanics through the effective force and peridynamic state concepts. The nonlocal deformation gradient is multiplicatively decomposed into the elastic and plastic parts in this large-deformation PPM paradigm. The stabilized correspondence principle is adopted to implement a classical elastoplastic constitutive model into the new PPM paradigm. We have numerically implemented this large-deformation plasticity PPM paradigm through a Lagrangian meshfree method in space and an explicit Newmark scheme in time. The implemented PPM framework is used to model shear banding and retrogressive slope failure in porous media under dry conditions. The numerical results have demonstrated the efficacy and robustness of this new PPM paradigm for modeling shear banding and retrogressive failure involving large deformation in porous media.

## 1. Introduction

Large deformation localized failure in geomaterials, such as shear banding and retrogressive failure, is ubiquitous in geohazards such as landslides [1-2]. For instance, retrogressive failure, in the context of landslides, refers to a type of landslide movement where the failure or collapse of the earth materials begins at the toe (lower portion) of the initial landslide and then progressively moves upslope, eroding and incorporating additional material as it advances. This process can result in the enlargement and extension of the landslide. Thus, the characterization and prediction of localized failure in geomaterials play a crucial role in geohazards engineering. In this article, we investigate shear bands and retrogressive failure of geomaterials through a recently formulated large-deformation-plasticity periporomechanics (PPM) paradigm, a nonlocal meshfree numerical method.

PPM is a nonlocal reformulation of classical poromechanics through the effective force and peridynamic state concepts [3-9]. The field equations of PPM are intego-differential equations (integration in space and differentiation in time) in lieu of partial differential equations. The classical constitutive models for geomaterials can be readily implemented in PPM through the stabilized correspondence principle [3]. Previously, following the small strain theory, the total strain is decomposed additively into an elastic part and a plastic one. In this newly formulated large-deformation-plasticity PPM paradigm, the nonlocal deformation gradient is multiplicatively decomposed into the elastic and plastic parts following the classical computational large-deformation plasticity [10]. We have numerically implemented this large-deformation-plasticity PPM paradigm through an explicit Lagrangian meshfree scheme.

In this article, we study the inception of shear banding in geomaterial under biaxial loading conditions and the retrogressive slope failure utilizing the large-deformation-plasticity PPM paradigm. Both numerical examples involve large deformation in geomaterials. In the example of shear banding, the nonlocal second-order work criterion is used to validate the inception of shear banding. In the second example, the numerical result of the retrogressive failure is compared with the data in the literature [11].

**2. Mathematical formulation**

This section presents the governing equations of the PPM paradigm for porous media under dry conditions (e.g., a single phase). In the PPM paradigm, the porous media is assumed to consist of a finite number of material points and the material points at a finite distance that is not greater than $\delta$ called horizon has direct mechanical interaction [1-5]. A mechanical bond between the material points at $x$ and $x'$ in the initial configuration is denoted as $\underline{\xi} = x' - x$. The equation of motion for a dry porous body is written as

$$\rho_s \ddot{u} = \int_{\mathcal{H}} (\underline{\overline{T}} - \underline{\overline{T}}') \, dv + \rho_s g \tag{1}$$

where $\rho_s$ is the partial density of solid phase, $\ddot{u}$ is the acceleration, $\mathcal{H}$ is a spherical domain with a radius of $\delta$ around the material point $x$, $\underline{T}$ and $\underline{T}'$ are the effective force state at material points $x$ and $x'$ along the bond $\underline{\xi}$, respectively, $g$ is the gravity acceleration. The partial density of the solid skeleton is defined as

$$\rho_s = (1 - \emptyset)\rho^s \tag{2}$$

where $\rho^s$ is the intrinsic density of the solid skeleton and $\emptyset$ is the porosity [4]. The deformation state $\underline{Y}$ and the displacement state $\underline{U}$ associated with the bond $\underline{\xi}$
are defined as

$$\underline{U} = u' - u \quad \text{and} \quad \underline{Y} = y' - y \tag{3}$$

where $u'$ and $u$ are the displacements of material points $x'$ and $x$, respectively, $y'$ and $y$ are the location of material points $x'$ and $x$, respectively, in the deformed configuration.

Given $\underline{Y}$, the deformation gradient tensor $F$ in PPM is defined as

$$F = \left[ \int_{\mathcal{H}} \underline{\omega} \left( \underline{Y} \otimes \underline{\xi} \right) dv \right] \mathcal{K}^{-1} \tag{4}$$

where $\underline{\omega}$ is the influence function and $\mathcal{K}$ is the shape tensor [2]. Through the correspondence principle, the effective force state in PPM can be written as

$$\underline{\overline{T}} = \underline{\omega} \, \overline{P} \mathcal{K}^{-1} \underline{\xi} \tag{5}$$

where $\overline{P}$ is the effective Piola stress tensor. It follows from (5) that the effective force state can be computed from the classical constitutive models given the deformation state.

*2.1. Kinematics of large deformation elastoplasticity*

For the finite deformation plasticity, a local multiplicative decomposition of the deformation gradient $F$ is written as

$$F = F^e F^p \tag{6}$$

where $F^e$ is the elastic part of the deformation gradient and $F^p$ is the plastic part of the deformation gradient. The total and elastic left Cauchy-Green tensors can be written as

$$b = F F^T \quad \text{and} \quad b^e = F^e F^{eT} \tag{7}$$

The plastic right Cauchy-Green tensor is written as

$$C^p = F^{pT} F^p \tag{8}$$

The rotational tensor is defined as

$$R = VF^{-1} \text{ and } V = \sqrt{b} \tag{9}$$

In the PPM framework, the time derivative of deformation gradient tensor is defined as

$$\dot{F} = \left[ \int_{\mathcal{H}} \omega \left( \underline{\dot{U}} \otimes \underline{\xi} \right) dv \right] \mathcal{K}^{-1} \tag{10}$$

From the time derivative of the deformation gradient tensor, the spatial velocity gradient is written as

$$L = \dot{F} F^{-1} \tag{11}$$

The spatial rate of the deformation tensor can be written as

$$d = \frac{1}{2}(L + L^T) \tag{12}$$

The strain increment can be obtained from the spatial rate of deformation tensor.–Figure 1 shows the kinematics for the large deformation PPM framework in the reference, current configuration, and next configuration for the two material points $x$ and $x'$. In Figure 1, $f$ is the relative deformation gradient tensor.

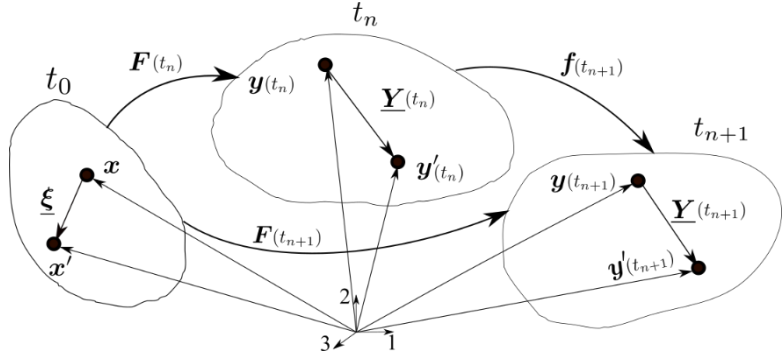

**Figure 1.** Kinematics of the material points $x$ and $x'$ in the reference, current configuration, and next configuration in the PPM

*2.2. Classical elastoplastic material model*

The effective Kirchhoff stress $\tau$ [10] can be obtained from the classical elastoplastic model. The corotational Kirchhoff stress tensor $\bar{\tau}$ is written as

$$\bar{\tau} = R^T \tau R \tag{13}$$

The Piola stress tensor can be written from the corotational Kirchhoff stress tensor as

$$\bar{P} = \bar{\tau} F^{-T} \tag{14}$$

The Lie derivative of the Kirchhoff stress tensor can be obtained from spatial rate of deformation tensor as

$$\mathcal{L}_v \bar{\tau} = a : d^e \text{ and } a = \lambda (\mathbf{1} \otimes \mathbf{1}) + 2G I \tag{15}$$

where $\lambda$ is the first Lame's constant, $G$ is the shear modulus, $\mathbf{1}$ is the second-order identity tensor, and $I$ is the fourth-order symmetric identity tensor.

*2.2.1 Drucker-Prager model*

In this study, the Drucker-Prager yield criteria is adopted, which is written as

$$\mathcal{F} = q + \alpha_1 \bar{p} - \alpha_2 c \tag{16}$$

where $c$ is the cohesion, $\alpha_1$ and $\alpha_2$ are material parameters that depend on the frictional angle, $q$ is the deviatoric stress, and $\bar{p}$ is the effective mean stress. The effective mean stress and deviatoric stress are defined as

$$\bar{p} = \frac{1}{3}(\bar{\tau}:\mathbf{1}) \text{ and } q = \sqrt{\frac{3}{2}}\|\mathbf{s}\| \qquad (17)$$

where $\mathbf{s}$ is the shear stress tensor. For the isotopic hardening $c$ is a function of the internal hardening variable $\zeta$ as

$$c = c_0 + h\zeta \qquad (18)$$

where $h$ is the isotropic hardening modulus, $c_0$ is the initial cohesion. The plastic potential is defined as

$$\mathcal{G} = q + \alpha_3 \bar{p} - \alpha_4 c \qquad (19)$$

where $\alpha_3$ and $\alpha_4$ are the material parameters that depend on the dilatancy angle.

*2.2.2 Flow rule and hardening law*
For non-associated plasticity, the evolution equation can be obtained from the plastic potential as

$$\dot{\mathbf{C}}^p = -2\dot{\gamma}\frac{\partial \mathcal{G}}{\partial \tau}\mathbf{F}^{-1}\mathbf{b}^e\mathbf{F}^{-T} \qquad (20)$$

The rate form of the left Cauchy-Green strain tensor can be written as

$$\dot{\mathbf{b}}^e = -2\dot{\gamma}\frac{\partial \mathcal{G}}{\partial \tau}\mathbf{b}^e \qquad (21)$$

For the hardening law, the rate of the internal hardening variable [10] can be written as

$$\dot{\zeta} = \dot{\gamma}\frac{\partial \mathcal{G}}{\partial Q} \qquad (22)$$

where $Q$ is a thermodynamic variable conjugated to the internal hardening variable $\zeta$. We refer to [10] for more technical details on the numerical implementation of a finite strain plasticity model.

*2.3 Numerical implementation*
The large deformation plasticity PPM paradigm has been implemented numerically through a Lagrangian meshfree method in space and an explicit Newmark scheme in time [3]. In the subsequent section, we present two numerical results on localized failures in dry porous media through the implemented large-deformation-plasticity PPM paradigm.

**3. Numerical examples**
This section presents two numerical examples of localized failure with large deformation in porous media. Example 1 deals with conjugate shear bands. Example 2 concerns the retrogressive slope failure.

*3.1. Example 1: Conjugate shear bands*
This example simulates the formation of conjugate shear bands in porous plastic media under dynamic loading conditions. Figure 2 plots the model setup for this example. A vertical displacement $u_y = 0.4$ m is applied on the top boundary with a rate of 0.4 m/s. The constant lateral confining pressure 0.1 MPa is applied on the left and right boundaries. The input material parameters are: solid phase density $\rho$ =2000 kg/m³, initial porosity $\phi_0$=0.2, bulk modulus $K$ =3.8 MPa, shear modulus $G$ =2.2 MPa, initial cohesion $c_0$= 20 kPa, residual cohesion $c_r$=8 kPa, softening modulus $h = -20$ kPa, frictional angle $35°$, and dilatation angle $15°$. The gravity load is not considered. The specimen domain is discretized into $40 \times 80$ material points with $\Delta x$=0.75 m, and horizon size is $\delta= 3\Delta x$. The time increment is $\Delta t = 1 \times 10^{-4}$ s.

The numerical results are presented in Figures 3 - 7. Figure 3 plots the loading curve on the top boundary. The results show that the softening occurs at $u_y = 0.1$ m. Figure 4 plots the snapshots of equivalent plastic shear strain in the deformed configurations at three loading stages. The results in Figure 4 show that the conjugate shear bands have developed gradually. Figure 5 presents the snapshots of plastic volumetric strain at the three loading stages. The positive sign of the plastic volume strain denotes dilatation. Figure 6 shows the snapshots of the displacement magnitude at the same three loading stages. The nonlocal second-order work criterion can be used to validate the numerical results of the localized instability, such as shear bands in porous media [8]. Figure 7 shows the snapshots of second-order work [8] in the deformed configurations of the three loading stages. The results in Figure 7 demonstrate that the zone of the negative second-order work is consistent with the shear banding zone in the specimen.

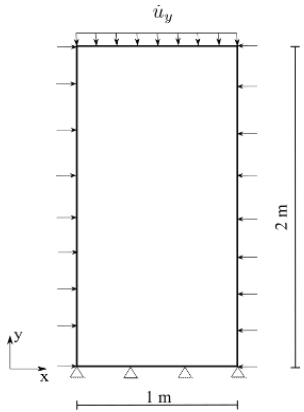

**Figure 2.** Model setup for the shear banding example.

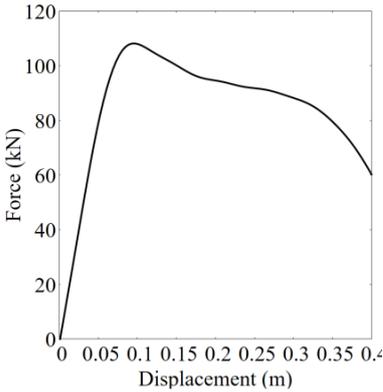

**Figure 3.** Loading curve on the top boundary.

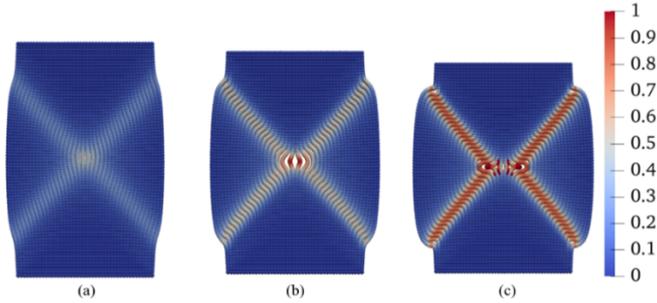

**Figure 4.** Contours of equivalent plastic shear strain in the deformed configurations at (a) $u_y = 0.1$, (b) $u_y = 0.25$, and (c) $u_y = 0.4$.

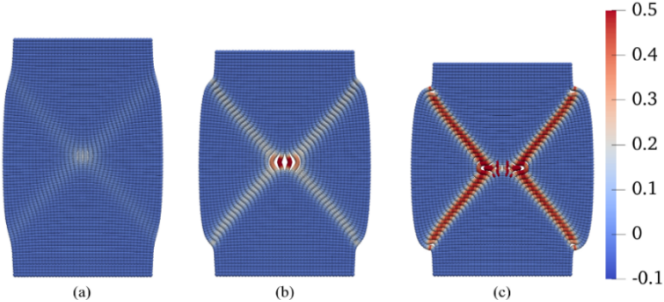

**Figure 5.** Contours of the plastic volume strain in the deformed configurations at (a) $u_y = 0.1$, (b) $u_y = 0.25$, and (c) $u_y = 0.4$.

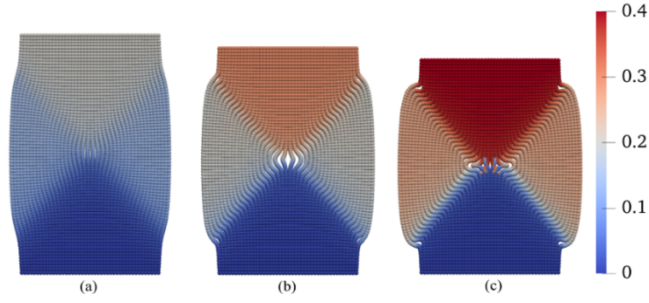

**Figure 6.** Contours of the displacement magnitude in the deformed configurations at (a) $u_y = 0.1$, (b) $u_y = 0.25$, and (c) $u_y = 0.4$.

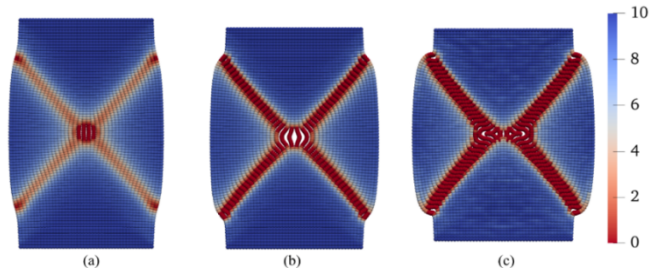

**Figure 7.** Contours of the second order work in the deformed configurations at (a) $u_y = 0.1$, (b) $u_y = 0.25$, and (c) $u_y = 0.4$.

*3.2. Example 2: Retrogressive slope failure*

This example simulates the retrogressive slope failure reported in [11]. In 1994, a landslide on sensitive clay occurred in Sainte-Monique, Quebec. The retrogression distance is about 100 m from the initial crest to the back scarp. The landslide is characterized as a spread failure with horsts and grabens observed in the site. The spread failures resulted from the extension and dislocation of the soil mass. Direct field measurements showed a horst top angle of about $52°$.

Figure 8 shows the problem setup of this example. The height of the slope is 17 m, and the inclination angle is $24°$. At the bottom boundary, the friction factor between the slope and the base is 0.3. The left side boundary is fixed in the horizontal direction. The eroded part is shown as dashed lines at the slope toe, as shown in Figure 8. The initial stress condition is generated by gravity loading with $K_0 = 0.5$. The failure is triggered by removing the retaining force on the slope surface. The material parameters adopted for this example are: solid phase density $\rho_s$=1600 kg/m³, initial porosity $\phi_0$=0.2, Young's modulus $E$ =1 MPa, Poisson's ratio $\nu$= 0.495, initial cohesion $c_0$=35 kPa, residual cohesion $c_r$=10 kPa, softening modulus $h = -1$ kPa, frictional angle $0°$, and dilatation angle $0°$. The domain is discretized into 10800 material points with $\Delta x$=0.75 m, and horizon size $\delta$= 3$\Delta x$. The time increment is $\Delta t = 1 \times 10^{-4}$.

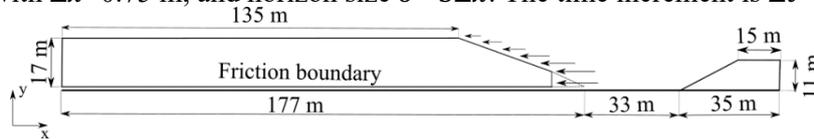

**Figure 8**. Model setup for the retrogressive slope failure example.

The numerical results are presented in Figures 9 and 10. Figure 9 plots the contours of the equivalent plastic shear strain in the deformed configurations at 5 s, 10 s, 17 s, 25 s, and 38 s, respectively. The results show that the shear band develops from the bottom boundary toward up. Five prominent collapses during the landslide can be observed, with the final configuration obtained at 38 s. As shown in Figure 9, the tip of the horsts in the simulation is about $60°$, which is comparable with the angle of horsts in the field. The difference can be related to the friction between the slope surface and the base.

Figure 10 shows the contours of the displacement magnitude in the deformed configurations at 5 s, 10 s, 17 s, 25 s, and 38 s, respectively. Figure 11 compares the ground levels in the site from the numerical results with the field data. In case 1, $c_0$= 30 kPa, and in case 2, $c_0$= 35 kPa, respectively. The other parameters remain the same. In this example, the large-deformation-plasticity PPM model is robust in simulating the formation and propagation of the inclined shear bands and reproduces the horsts and grabens that represent the spread failure mode.

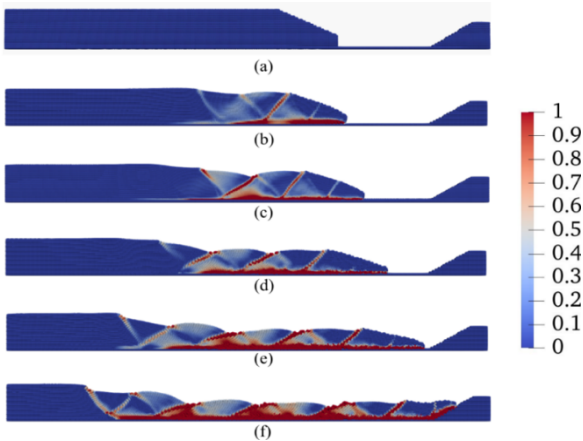

**Figure 9**. Contours of the equivalent plastic shear strain in the deformed configurations at times: (a) 0 s, (b) 5 s, (c) 11 s,(d) 17 s, (e) 25 s, and (f) 38 s.

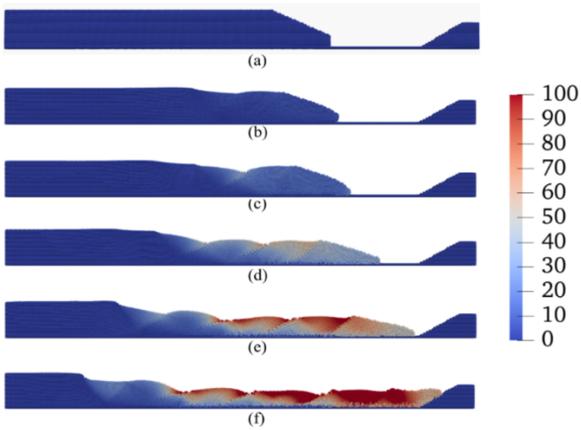

**Figure 10.** Contours of the displacement magnitude (m) in the deformed configurations at times: (a) 0 s, (b) 5 s, (c) 11 s, (d) 17 s, (e) 25 s, and (f) 38 s.

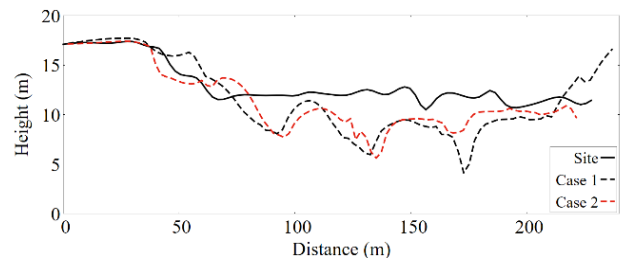

**Figure 11.** Comparison of ground levels for the Sainte-Monique landslide between the measurement and the simulation results.

## 4. Conclusion

This article presents a large-deformation-plasticity PPM paradigm for modeling shear banding and retrogressive failure in geomaterials. The classical elastoplasticity model is implemented into the framework through the multiplicative decomposition of the nonlocal deformation gradient following the computational finite strain plasticity theory. The large-deformation-plasticity PPM paradigm has been numerically implemented through an explicit Lagrangian mesh-free algorithm. The numerical examples have demonstrated the robustness and efficacy of the proposed PPM paradigm for modeling large-deformation failure in geomaterials. For instance, the new numerical paradigm can realistically capture the retrogression and run-out distance characteristics of the slope failure that manifests shear bands and large deformation.

**Acknowledgments**

Support for this work was provided by the US National Science Foundation under contract number CMMI 1944009.